\documentclass{article}
\usepackage{graphicx}

\usepackage[fleqn]{amsmath}
\usepackage{amsthm}
\usepackage{amsmath}
\usepackage{amssymb}
\usepackage{amsbsy}
\usepackage{amsfonts}
\usepackage{mathrsfs}
\usepackage[all]{xy}
\usepackage{amstext}
\usepackage{amscd}
\usepackage[dvips]{epsfig,color}
\usepackage{psfrag}
\usepackage{enumerate}
\usepackage{flafter}
\allowdisplaybreaks

\theoremstyle{plain}
\newtheorem{thm}{Theorem}[section]
\newtheorem{prop}[thm]{Proposition}
\newtheorem{cor}[thm]{Corollary}

\theoremstyle{definition}
\newtheorem{exa}[thm]{Example}

\newtheorem{rem}[thm]{Remark}

\newtheorem{defbn}[thm]{Definition}

\def\Ker{\mathop{\mathrm{Ker}}\nolimits}

\def\Hom{\mathop{\mathrm{Hom}}\nolimits}

\newcommand{\lra}{\longrightarrow}
\newcommand{\ra}{\rightarrow}
\newcommand{\Q}{{\Bbb Q}}
\newcommand{\R}{{\Bbb R}}
\newcommand{\Z}{{\Bbb Z}}
\newcommand{\N}{{\Bbb N}}

\newcommand{\PD}{\mathrm{PD} }

\newcommand{\DW}{\mathrm{DW}}

\newcommand{\pc}[2]{\mbox{$\begin{array}{c}
\includegraphics[scale=#2]{#1.pdf}
\end{array}$}}

\begin{document}

\large
\begin{center}{\bf\Large Cellular chain complexes of universal covers of some 3-manifolds}
\end{center}
\vskip 1.5pc
\begin{center}{Takefumi Nosaka\footnote{
E-mail address: {\tt nosaka@math.titech.ac.jp}
}}\end{center}
\vskip 1pc
\begin{abstract}\baselineskip=12pt \noindent
For a closed 3-manifold $M$ in a certain class, we give a presentation of the cellular chain complex of the universal cover of $M$. The class includes all surface bundles, some surgeries of knots in $S^3$, some cyclic branched cover of $S^3$, and some Seifert manifolds. In application, we establish a formula for calculating the linking form of a cyclic branched cover of $S^3$, and develop procedures of computing some Dijkgraaf-Witten invariants.
\end{abstract}
\begin{center}
\normalsize
\baselineskip=17pt
{\bf Keywords} \\
\ \ \ Universal covering, 3-manifold group, group homology, knot, branched coverings, linking form \ \
\end{center}

\large
\baselineskip=16pt
\section{Introduction}
\label{IntroS}
In order to investigate a connected CW-complex $X$ with a non-trivial fundamental group $\pi_1(X)$, it is important to give a concrete presentation of the cellular chain complex, $C_*(\widetilde{X};\Z)$, and the cup-products of the universal cover $\widetilde{X}$. In fact, the homology of $X$ with local coefficients and the (twisted) Reidemeister torsion of $X$ are defined from $C_*(\widetilde{X};\Z)$. If $X$ is a $K(\pi,1)$-space, the chain complex means a projective resolution of the group ring $\Z[\pi_1(X)] $.
Thus, it is also of use for computing many invariants to concretely present $ C_*(\widetilde{X};\Z). $

This paper focuses on a class of closed 3-manifolds satisfying the following condition:

\vskip 0.517pc
\noindent
{\bf Assumption ($\dagger$)} 
A closed oriented 3-manifold $M$ satisfies that any closed 3-manifold $M'$ with a group isomorphism $\pi_1(M) \cong \pi_1(M')$ admits a homotopy equivalence $M \simeq M'$.

\vskip 0.517pc

\noindent For example, $M$ satisfies this assumption if $M $ is an Eilenberg-MacLane space of type $(\pi_1(M),1)$, which is equivalent to that $M$ is irreducible and has an infinite fundamental group. In Section \ref{defS3}, we examine many 3-manifolds, including all surface bundles, some surgeries of knots in $S^3$, spliced sums, cyclic branched covers of $S^3$ with Assumption ($\dagger$), and some Seifert manifolds. For when $M$ is one of these, we describe presentations of the complex $C_*(\widetilde{M};\Z)$ and of the cup-product $H^1(M;N)\otimes H^2(M;N') \ra H^3(M;N \otimes N') $ for any local coefficient modules $N,N'$. The procedure for obtaining such descriptions essentially follows from the work of \cite{Sie1,Tro} in terms of ``identity", which we review in Section \ref{ganeaS}. This procedure can also be used to describe the fundamental homology 3-class, $[M]$ of $M$; see Remark \ref{bb22}.

In application, we give a formula for the linking forms of cyclic branched covers of $S^3$ with Assumption ($\dagger$) (see Propositions \ref{bbbn1}).
Furthermore, we develop procedures of computing some Dijkgraaf-Witten invariants from 
the above descriptions; see \S \ref{defS42}.
In addition, such descriptions of identities are used for computing knot concordance groups, Reidemeister torsions, and Casson invariants; see \cite{MP,Nos,W}. 
There might be other applications from the above presentations of the complexes $C_*(\widetilde{M};\Z)$

\noindent
{\textbf{Conventional notation.} In this paper, every manifold is understood to be smooth, connected, and orientable. By $M $, we mean a closed 3-manifold with orientation $[M]$.

\section{Taut identities and cup-products }
\label{ganeaS}
\subsection{Review: identities and cup-products }
\label{ganeaS1}
Let us recall the procedure of obtaining cellular chain complexes of some universal covers, as described in the papers \cite{Sie1} and \cite{Tro}.
There is nothing new in this section. 

We will start by reviewing identities. Take a finitely presented group $\langle x_1, \dots, x_m \mid r_1, \dots, r_m \rangle $ of deficiency zero. Setting up the free groups $F:=\langle x_1, \dots, x_m \mid \rangle $ and $P:=\langle \rho_1, \dots, \rho_m | \rangle $, let us consider the homomorphism,
\[ \psi: P* F \lra F \ \ \ \ \mathrm{defined \ by \ \ \ \ } \psi(\rho_j)=r_j, \ \ \ \psi (x_i)=x_i. \]
An element $s \in P*F $ is {\it an identity} if $s \in \Ker(\psi)$ and $s$ can be written as $ \prod_{k=1}^n \omega_k \rho_{j_k}^{\epsilon_k} \omega_k^{-1}$ for some $ w_k \in F$, $\epsilon_k \in \{ \pm 1\} $ and indices ${j_k}$'s.

Given a closed 3-manifold $M$ with a genus-$m$ Heegaard splitting, let us review the cellular complex of the universal cover, $\widetilde{M}$, of $M.$ A CW-complex structure of $M$ induced by the splitting consists of a single zero-cell, $ m$ one-handles, $m$ two-handles, and a single three-handle. Therefore, $\pi_1(M)$ has a group presentation $\langle x_1, \dots, x_m \ \mid \ r_1, \dots, r_m \rangle $, and the cellular complex of $\widetilde{M}$ is described as
\begin{equation}\label{b78} C_*(\widetilde{M};\Z): 0 \ra \Z[\pi_1(M)] \stackrel{\partial_3}{\lra} \Z[\pi_1(M)]^m \xrightarrow{ \ \partial_2 \ } \Z[\pi_1(M)]^m \stackrel{\partial_1 }{\lra} \Z[\pi_1(M)] \ra 0.\end{equation}
Here, $\Z[\pi_1(M)] $ is the group ring of $\pi_1(M). $ We will explain the boundary maps $\partial_*$ in detail. Let $\{ a_1, \dots, a_m\}$, $\{ b_1, \dots, b_m\} $, and $\{ c\}$ denote the canonical bases of $C_1(\widetilde{M};\Z)$, $C_2(\widetilde{M};\Z)$, and $C_3(\widetilde{M};\Z)$ as left $\Z[\pi_1(M)]$-modules, respectively. Then, as is shown in \cite{Lyn}, $\partial_1(a_i)=1-x_i$, and $\partial_2(b_i)=\sum_{k=1}^m [ \frac{\partial r_i}{\partial x_k}] a_k$, where $\frac{\partial r_i}{\partial x_{k}} $ is the Fox derivative. Moreover, the main result in \cite{Sie1} is that there exists an identity $s$ such that $ \partial_3(c)= \sum_k [\psi ( \frac{\partial s}{\partial \rho_k})] b_k$.

Next, we will briefly give a formula for the cup-product in terms of the identity, which is a result of \cite[\S 2.4]{Tro}. Let $N$ and $N'$ be left $\Z[ \pi_1(M)]$-modules. We can define the cochain complex on $C^*( M;N):= \Hom_{\Z[ \pi_1(M)] }(C_* (\widetilde{M};\Z) , N ) $ with local coefficients. Recalling the definition of the identity $s=\prod_{k=1}^n \omega_k \rho_{j_k}^{\epsilon_k} \omega_k^{-1}$, define
\begin{equation}\label{b4} D^{\sharp}(c) = \sum_{k=1}^n \epsilon_k \bigl( \sum_{\ell=1}^m [ \frac{\partial \omega_k}{\partial x_\ell}] a_\ell \otimes \omega_k b_{j_k} \bigr) \in C_1(\widetilde{M};\Z) \otimes C_2(\widetilde{M};\Z) . \end{equation}
Then, for cochains $p \in C^1 ( M;N)$ and $q \in C^2 ( M;N')$, we define a 3-cochain $ p \smile q$ by
$$ p \smile q (u c):= (p \otimes q)(u D^{\sharp} (c)) \in N \otimes_{\Z}N'. $$
Here, $u \in \Z[\pi_1(M)]$. Then, the map
$$\smile : C^1(M ;N) \otimes C^2(M;N') \ra C^3(M ;N \otimes_{\Z}N' ); \ \ (p,q) \mapsto p \smile q, $$
induces the bilinear map on cohomology, which is known to be equal to the usual cup-product. Here, notice that, since the third $\partial_3 \otimes_{\Z [\pi_1(M) ]} \mathrm{id}_{\Z} $ is zero, the 3-class $s \otimes 1 \in C_*(\widetilde{M}) \otimes_{\Z[\pi_1(M) ]} \Z $ is a generator of $H_3( C_*(\widetilde{M}) \otimes \Z) \cong H_3( M ; \Z) \cong \Z$, which represents the fundamental 3-class $[M]$; thus, given a $\pi_1(M)$-invariant bilinear map $\psi : N \otimes N' \ra A$ for some abelian group $A$, we have the following equality on the pairing of $[M]$:
\begin{equation}\label{b8w} \psi \circ \smile ( p, q) = \psi\langle p \smile q , [M]\rangle \in A, \end{equation}
for any cochains $p \in C^1 ( M;N)$ and $q \in C^2 ( M;N')$.

In summary, for a description of the complex $C_*(\widetilde{M})$ and the cup-product, it is important to describe an identity from $M$. 
\subsection{Taut identities}
\label{ganeaS2}
In order to find such identities giving the complex \eqref{b78}, we review tautness from \cite{Sie1}; see also \cite[Appendix]{W} for a brief explanation.
Fix a finite presentation $ \langle x_1, \dots, x_m \mid r_1, \dots, r_m \rangle $. 
Let $s= \prod_{k=1}^{2m} w_k \rho_{j_k}^{\epsilon_k} w_k^{-1} \in P * F $ be an identity, where
$\rho_{j_k}$ and $w_{k}$ can be written in
$$\rho_{j_k}= a_{k, 1}^{\epsilon_{k,1}} \cdots a_{k,\ell_k}^{\epsilon_{k ,\ell_k }} , \ \ \ \ \ \ w_{k} =b_{k,1}^{\eta_{k,1}} \cdots b_{k,n_k}^{\eta_{k, n_k}}, \ \ \ \ \ \ \ \ \ (\epsilon_{i,j}, \eta_{i,j} \in \{ \pm 1\}) .$$
Here, $a_{k, \ell}$ and $b_{k, \ell} $ lie in $\{ x_1, \dots, x_m \}$.
For each $w_k \rho_{j_k}^{\epsilon_k} w_k^{-1} $, take the $\ell_k$-gon
$D_{j_k}$ whose $i$-th edge is labeled by $a_{k,i}^{\epsilon_{k ,i }}$, and
the segment $I_k=[0,n_k ]$ such that $[ i-1,i]$ is labeled by $b_{k,i} ^{\eta_{k,i}}$.
\begin{defbn}[\cite{Sie1}]\label{bdays}
\begin{enumerate}[(1)]
\item A self-bijection$$\mathcal{I}: \cup_{k=1}^{2m} \{ (k,1),\dots, (k, \ell_k) \} \ra \cup_{k=1}^{2m} \{ (k,1),\dots, (k, \ell_k) \} $$
is called a {\it syllable} if $a_{\mathcal{I}(i,j) } = a_{i,j} \in F$ and $ \epsilon_{i,j} =- \epsilon_{\mathcal{I}(i,j)} \in \{ \pm 1\}$.

\item For a syllable $\mathcal{I}$, consider the following equivalence on the disjoint union $\sqcup_{ i=1}^{2m} D_{r_i }$: the interval with labeling $a_{i,j} $ is identified with those with labeling $a_{\mathcal{I}(i,j)} $.

\item An identity $s$ is said to be {\it taut} if there is a syllable $\mathcal{I}$ such that the quotient space $\sqcup_{ i=1}^{2m} D_{r_i } / \sim $ of $\sqcup_{ i=1}^{2m} D_{r_i }$ subject to the above equivalence $\sim $ is homeomorphic to $S^2 $, and if there are injective continuous maps
$$ \kappa_k : I_k=[0, n_{k} ] \ra \sqcup_{ i=1}^{2m} \partial D_{r_i } / \sim \textrm{,} \ \ \ \ \ \ \lambda_k : [0, \ell_{k} ) \ra \partial D_{r_k } /\sim $$
satisfying the following condition (*).

(*) For each $k$, the image $\kappa_k([i-1, i]) $ coincides with an edge labeled by $b_{k,i}$ compatible with the orientations, and $ \lambda_k ([ j-1, j]) $ coincides with the $j$-th edge of $ D_{r_k } $ compatible with the orientations.
Furthermore, $\kappa_k(n_k )= \lambda_k (0)= \lambda_k (\ell_k)$.
\end{enumerate}
\end{defbn}
This paper is mainly based on the following theorem of Sieradski: 
\begin{thm}[\cite{Sie1}]\label{bb2}
Given a group presentation $\langle x_1, \dots, x_m \ \mid \ r_1, \dots, r_m \rangle $ with a taut identity $s$, there exists a closed 3-manifold $M$ with a genus-$m$ Heegaard splitting such that the complex $C_*(\widetilde{M};\Z)$ is isomorphic to the complex \eqref{b78}. 
\end{thm}

In a concrete situation where an identity $s$ is explicitly described, it is not so hard to find such a $\mathcal{I}$ and show the tautness of $s$ (in fact, this check is to construct a 2-sphere from the disjoint union $\sqcup_{ i=1}^n D_{r_i }$ as a naive pasting).
In all the statements in \S \ref{defS3}, we will claim that some identities satisfy the taut condition; however, we will also omit the check by elementary complexity, as in other papers on taut identities \cite{BH,Sie1,Tro}.
\begin{exa}\label{bb244}
As an easy example of the pasting, we focus on the 3-dimensional torus $M=(S^1)^3$ with presentation $\pi_1(M)= \langle x,y,z \mid r ,s,u \rangle $,
where $ r =[x,y], s=[y,z],u=[z,x].$ As in \cite{Sie1}, consider the following identity.
$$W_{(S^1)^3}=r ( y^{-1} u^{-1} y) s ( z^{-1} r z) u (x^{-1 } s^{-1} x). $$
Then, Figure \ref{baki4} gives a self-bijection and $\lambda_m, \kappa_m $ satisfy the tautness. 
Moreover, if we attach a 3-ball in the right hand side in the figure along the boundary of the 3-cube, the resulting space is equal to $(S^1)^3$.
\end{exa}

\begin{figure}[h]
\begin{center}
\begin{picture}(100,110)


\put(-110,42){\pc{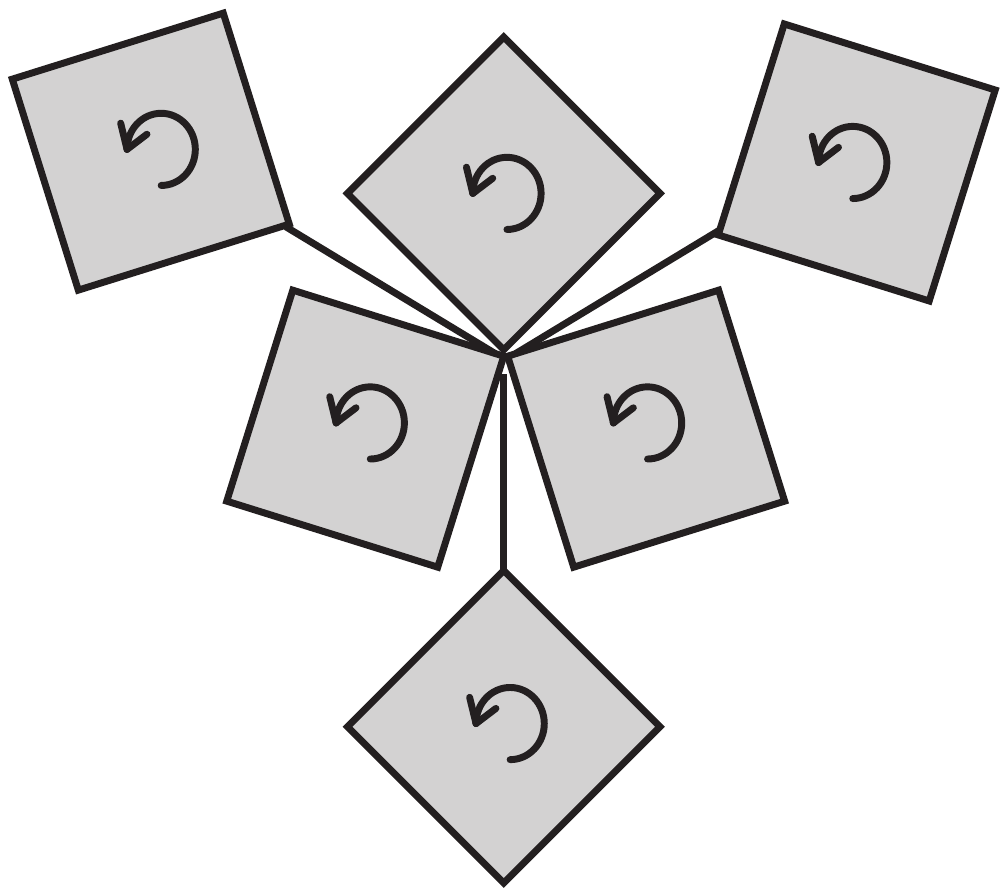}{0.36156}}

\put(95,-11){\pc{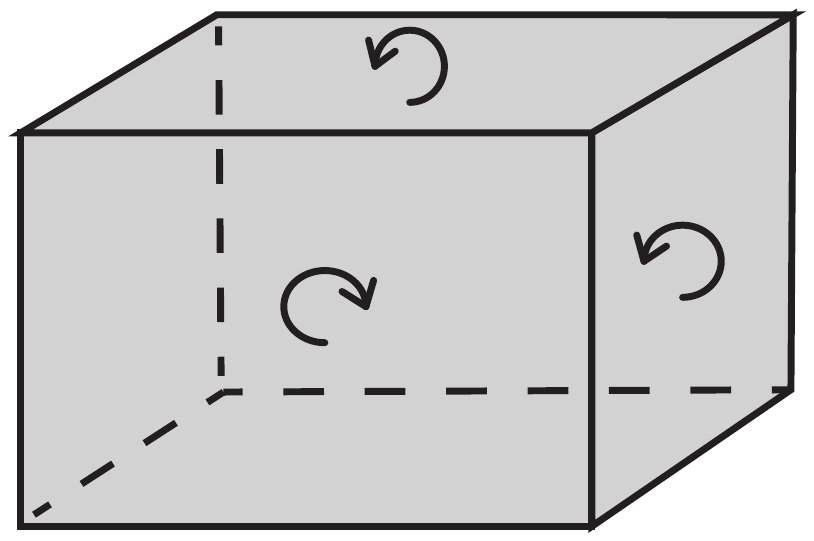}{0.4115156}}

\put(45,51){\large $\Phi$}
\put(35,44){\large $- \!\!\! - \!\!\!- \!\!\!\lra$}
\put(-31,54){ $D_1$}
\put(-96,54){ $D_5$}
\put(-6,74){ $D_2$}
\put(-121,74){ $D_4$}
\put(-76,84){ $D_3$}
\put(-40,68){ $I_2$}
\put(-79,68){ $I_4$}

\put(-102,33){ $I_6 - \!\!\! - \!\!\! - \!\!\!\!\! - \!\!\! - \!\!\!\!\!\lra $}

\put(-81,4){ $D_6$}

\end{picture}
\caption{The tautness of $ (S^1)^3$. The right side means the 2-sphere obtained as the quotient $\sqcup_{ i=1}^6 D_{i } / \sim $. Here,
the restriction map on $I_i $ of $\Phi$ means $\lambda_i$, and the restriction map on $\partial D_i $ of $\Phi$ means $\kappa_i$.}\label{baki4}
\end{center}
\end{figure}

\begin{rem}\label{bb22}
Suppose that we find a taut identity $s$ from $\langle x_1, \dots, x_m \ \mid \ r_1, \dots, r_m \rangle $, and 
the resulting 3-manifold $M$ satisfies Assumption ($\dagger$). Then, by Assumption ($\dagger$), the resulting 3-manifold up to homotopy does not depend on the choice of $s$. 
In particular, we emphasize that, if $M$ satisfies Assumption ($\dagger$) and we find a taut identity from $\pi_1(M)= \langle x_1, \dots, x_m \ \mid \ r_1, \dots, r_m \rangle $, then the third $\partial_3$ and the cup-product are uniquely determined, up to homotopy, by the identity. In fact, if we have another identity $\omega'$ and consider the associated $ C_*(\widetilde{M})'$, Assumption ($\dagger$) ensures a chain map $C_*(\widetilde{M}) \ra C_*(\widetilde{M})'$, which induces a homotopy equivalence.
\end{rem}

\section{Descriptions of taut identities of various 3-manifolds}
\label{defS3}
In this section, we give several examples of identities from some classes of 3-manifolds. We will describe the cellular complexes of some universal covers.

\subsection{Fibered 3-manifolds with surface fibers over the circle}
\label{defS31}
First, we will focus on surface bundles over $S^1$. Let $\Sigma_{g}$ be an oriented closed surface of genus $g$ and $f: \Sigma_{g} \ra \Sigma_{g}$ an orientation-preserving diffeomorphism. {\it The mapping torus}, $T_{f}$, is the quotient space of $\Sigma_{g} \times [0,1]$ subject to the relation $(y,0) \sim (f(y),1)$ for any $y \in \Sigma_{g}$. The homeomorphism type of $T_f$ depends on the mapping class of $f$. Conversely, if a closed 3-manifold $M$ is a fibered space over $S^1$, then $M$ is homeomorphic to $T_{f}$ for some $f$. Since $T_{f} $ is a $\Sigma_{g}$-bundle over $S^1$, it is a $K(\pi,1)$-space and therefore satisfies Assumption ($\dagger$).

We will construct an identity. Choose a generating set $\{ x_1,\dots, x_{2g} \} $ of $\pi_1(\Sigma_{g})$, which gives the isomorphism $\pi_1(\Sigma_{g}) \cong \langle x_1,\dots, x_{2g} \mid [x_1, x_2] \cdots [x_{2g-1},x_{2g}]\ \rangle . $ Following a van Kampen argument, we can verify the presentation of $ \pi_1(T_{f})$ as
\begin{equation}\label{m3as2} \langle \ x_1, \dots, x_{2g} , \gamma \mid r_i:= \gamma f_*(x_i) \gamma^{-1} x_i^{-1}, \ \ \ ( i \leq 2g), \ \ r_{2g+1}:=[x_1, x_2] \cdots [x_{2g-1},x_{2g}] \rangle .
\end{equation}
Here, $\gamma$ represents a generator of $\pi_1(S^1). $ For $ i \leq 2g $, define $w_{i}= \prod_{j=1}^i [ x_{2j-1}, x_{2j} ]\in F $, and
\[ W_{i} := w_{i-1} \rho_{2i-1} w_{i-1}^{-1} \cdot
(w_{i-1} x_{2i-1}) \rho_{2i }(w_{i-1} x_{2i-1})^{-1}
\cdot (w_{i} x_{2i}) \rho_{2i-1}^{-1}(w_{i} x_{2i})^{-1} \cdot w_{ i} \rho_{2i}^{-1}w_{i}^{-1} . \]
Since $f$ can be isotoped so as to preserve a point $z\in \Sigma_g$, we regard the induced map $f_*$ as
a homomorphism $: \pi_1(\Sigma_{g} \setminus \{ z\} ) \ra \pi_1(\Sigma_{g} \setminus \{ z\})$.
Since $ f_*$ is a group isomorphism, there exists a unique element $q_f \in \langle x_1,\dots, x_{2g} | \rangle $ satisfying
$$ f_* ([x_1, x_2] \cdots [x_{2g-1},x_{2g}])=q_f ([x_1, x_2] \cdots [x_{2g-1},x_{2g}]) q_f^{-1} \in \langle x_1,\dots, x_{2g} | \rangle .$$
\begin{thm}\label{b3378}
Let $W$ be $ (\Pi_{i=1}^g W_{i}) \rho_{2g+1} (\gamma q_f \rho_{2g+1 }^{-1}q_f^{-1} \gamma^{-1}) \in F*P$. Then, $W$ is an identity.
\end{thm}
\begin{proof}
Direct calculation gives $ \psi( W_{i})= w_{i-1}\gamma [f_*(x_{2i-1}), f_*(x_{2i}) ] \gamma ^{-1}w_{i}^{-1}$, which implies
$$ \psi ( \Pi_{i=1}^g W_{i})= \gamma ( \Pi_{i=1}^g[ f_*(x_{2i-1}), f_*(x_{2i}) ] )\gamma ^{-1}w_{g}^{-1} =
\gamma \Pi_{i=1}^g[ f_*(x_{2i-1}), f_*(x_{2i}) ]\gamma ^{-1} ( \Pi_{i=1}^g[x_{2i-1}, x_{2i } ])^{-1} . $$
Hence, $\psi(W)=1$ by definition; that is, $W$ turns out to be an identity.
\end{proof}
\noindent
Furthermore, we can verify that $W$ is taut by the definition of $W$. Hence, from the discussion in \S \ref{ganeaS}, we can readily prove the following corollary.
\begin{cor}\label{co1}
Under the above terminology, the cellular chain complex of $\widetilde{T_f}$ is given by
\[C_*(\widetilde{T_f};\Z): 0 \ra \Z[\pi_1(T_f)] \stackrel{\partial_3}{\lra} \Z[\pi_1(T_f)]^{2g+1} \xrightarrow{ \ \partial_2 \ } \Z[\pi_1(T_f)]^{2g+1} \stackrel{\partial_1 }{\lra} \Z[\pi_1(T_f)] \ra 0.\]
Here, $\partial_1(a_i)=1-x_i$, $\partial_1(\gamma )=1-\gamma $, and $\partial_2$ and $\partial_3 $ have the matrix presentations,
\[ \left(
\begin{array}{cc}
\Bigl\{ \gamma \frac{\partial f_*(x_i)}{\partial x_{j}} -\delta_{ij} \Bigr\}_{1 \leq i,j \leq 2g} & \Bigl\{ 1- x_i \Bigr\}_{1 \leq i \leq 2g}^{\rm transpose} \\
\Bigl\{ \frac{\partial r_{2g+1}}{\partial x_{j}} \Bigr\}_{1 \leq j \leq 2g} & 0 \\
\end{array}
\right) , \]
\[ \left( \begin{array}{cc}
\Bigl\{ w_{j-1} -w_j x_{2j} , w_{j-1} x_{2j-1}- w_{j} \Bigr\}_{1 \leq j \leq g}, & 1 - \gamma q_f \\
\end{array}
\right). \]
Furthermore, the diagonal map $D^{\sharp}(c)$ is represented by
\[\Bigl( \sum_{i=1}^{g}\sum_{k=1}^{2g}\frac{\partial w_{i-1} }{\partial x_k}a_k \otimes w_{i-1}b_{2i-1}- \frac{\partial (w_{i}x_{2i-1} )}{\partial x_k}a_k \otimes w_{i}x_{2i} b_{2i-1}- \frac{\partial (w_{i-1}x_{2i -1 }) }{\partial x_k}a_k \otimes w_{-1i}x_{2i-1}b_{2i} \]
\[ \ \ \ \ + \frac{\partial w_{i}}{\partial x_k} a_k \otimes w_i b_{2i} \Bigr)- \bigl( \sum_{k=1}^{2g} \frac{\partial (\gamma q_f )}{\partial x_k}a_k \otimes \gamma q_f b_{2g+1} \bigr)+ a_{2g+1} \otimes (1- \gamma q_f) b_{2g+1}. \]
\end{cor}
\begin{rem}\label{re7438}
Corollary \ref{co1} for every $g$ is a generalization of the result of \cite{Mar}; the paper gives the cellular complexes of $\widetilde{T_f}$ only in the case $g=1$. We can verify that Corollary \ref{co1} with $g=1$ coincides with the results in \cite{Mar}.
\end{rem}
Finally, we mention the virtually fibered conjecture, which was eventually proven by Wise; see, e.g., \cite{BW}. This conjecture states that every closed, irreducible, atoroidal 3-manifold $M$ with an infinite fundamental group has a finite cover, which is homeomorphic to $T_f$ for some $f$. Let $d \in \N$ be the degree of the covering. Then, if we can find such a cover $p: T_f \ra M$, the pushforward of the above identity $W$ gives an algebraic presentation of $d [M]$.

\subsection{Spliced sums and $(p/1)$- and $(1/q)$-surgeries of $S^3$ along knots}
\label{defS33}
We will focus on spliced sums and some surgeries of $S^3$ along knots and construct taut identities. This section supposes that the reader has basic knowledge of knot theory, as in \cite[Chapters 1--11]{Lic}.

Let us review spliced sums. Take two knots $K, K' \subset S^3$ and an orientation-reversing homeomorphism $h : \partial (S^3 \setminus \nu K) \ra \partial (S^3 \setminus \nu K') $, where $\nu K$ means an open tubular neighborhood of $K$. Then, we can define a closed 3-manifold, $ \Sigma_h (K,K')$, as the attaching space $(S^3 \setminus \nu K) \cup_{h} (S^3 \setminus \nu K') $ with $ \partial (S^3 \setminus \nu K)$ glued to $\partial (S^3 \setminus \nu K') $ by $h$. This space is commonly referred to as {\it the spliced sum of $(K,K')$ via $h$}. Spliced sums sometimes appear in discussions on additivity of topological invariants; see, e.g., \cite{BC}. Further, choose the preferred meridian-longitude pair $(\mathfrak{m}, \mathfrak{l})$ (resp. $(\mathfrak{m}', \mathfrak{l}'$)) as a generating set of $\pi_1 \partial (S^3 \setminus \nu K) $ (resp. of $\partial (S^3 \setminus \nu K') $). If $ h_* : \pi_1 \partial (S^3 \setminus \nu K) \ra \pi_1 \partial (S^3 \setminus \nu K')$ is represented by $ \begin{pmatrix} 0 & 1 \\ 1 & p \\ \end{pmatrix} $ (resp. $ \begin{pmatrix} 1 & 0 \\ q & -1 \\ \end{pmatrix} $) for some $ p,q \in \Z$, we denote $ \Sigma_h (K,K') $ by $ \Sigma_{p/1} (K,K')$ (resp. $ \Sigma_{1/q} (K,K')$). In particular, if $K'$ is the unknot, then $ \Sigma_{p/1} (K,K')$ and $ \Sigma_{1/q} (K,K')$ are the closed 3-manifolds obtained by $(p/1)$- and $(1/q)$-Dehn surgery on $K$ in $S^3$, respectively.

Since the identities of $ \Sigma_{p/1} (K,K')$ and $\Sigma_{1/q} (K,K') $ will be constructed in an analogous way to \cite[Page 481]{Tro}, let us review the terminology in \cite{Tro}. Choose a Seifert surface $\Sigma$ of genus $g$ and a bouquet of circles $W \subset \Sigma$ such that $W$ is a deformation retract of $\Sigma$ and $\pi_1 (S^3 \setminus \Sigma)$ is a free group. 
21. Page}. For example, any Seifert surface obtained by a Seifert algorithm admits such a bouquet. Choose a bicollar $\Sigma \times [-1,1]$ of $\Sigma$ such that $ \Sigma \times \{ 0\} =\Sigma$. Let $\iota_{\pm}: \Sigma \ra S^3 \setminus \Sigma$ be embeddings whose images are $ \Sigma \times \{ \pm 1\}$. Take generating sets $\{ v_1, \dots, v_{2g} \}$ of $\pi_1 \Sigma $ and $ \{ x_1, \dots, x_{2g} \} $ of $ \pi_1(S^3 \setminus \Sigma)$, and set $u_i^{\sharp}:=(\iota_+)_* (v_i)$ and $u_i^{\flat}:=(\iota_-)_* (v_i)$, where we may suppose that $[v_1,v_2] \cdots [v_{2g-1},v_{2g} ]$ represents a loop of $\pi_1\partial \Sigma $; a van Kampen argument yields a presentation
\begin{equation}\label{oo4656} \langle \ x_1, \dots, x_{2g}, \mathfrak{m} \mid r_i:=\mathfrak{m} u_i^{\sharp} \mathfrak{m}^{-1} (u_i^{\flat})^{-1} \ \ \ (1 \leq i \leq 2g) \ \ \rangle \end{equation}
of $ \pi_1(S^3 \setminus K)$. Here, $\mathfrak{m} $ is a representative of a meridian in $ \pi_1(S^3 \setminus K)$, and the $x_i$'s lie in the commutator subgroup of $ \pi_1(S^3 \setminus K)$. Since the boundary loops of $ \pi_1 \Sigma $ and $\pi_1(S^3 \setminus \Sigma) $ are equal by definition, we should notice
\begin{equation}\label{b652} [ u_1^{\flat}, u_2^{\flat} ] \cdots [ u_{2g-1}^{\flat}, u_{2g}^{\flat} ]
= [ u_1^{\sharp}, u_2^{\sharp} ] \cdots [ u_{2g-1}^{\sharp}, u_{2g}^{\sharp} ] \in \pi_1(S^3 \setminus \Sigma) , \end{equation}
which we denote by $\mathfrak{l}$. In other words, $\mathfrak{l}$ means a preferred longitude of $K.$

In a parallel way, concerning the other $K'$, we have a generating set $ \{ x_1', \dots, x_{2g'}' \} $ of $ \pi_1(S^3 \setminus \Sigma')$ and can define appropriate words $\mathfrak{u}_i^{\sharp}$ and $\mathfrak{u}_i^{\flat}$ such that
\begin{equation}\label{oo463356} \pi_1(S^3 \setminus K') \cong \langle \ x_1', \dots, x_{2g'}', \mathfrak{m}' \mid r_i':= \mathfrak{m}' \mathfrak{u}_i^{\sharp} (\mathfrak{m}')^{-1} (\mathfrak{u}_i^{\flat})^{-1} \ \ \ (1 \leq i \leq 2g') \ \ \rangle. \notag \end{equation}
We also redefine $\mathfrak{l}'$ by $[ \mathfrak{u}_1^{\flat}, \mathfrak{u}_2^{\flat} ] \cdots [ \mathfrak{u}_{2g-1}^{\flat}, \mathfrak{u}_{2g}^{\flat} ]. $

Before we state Theorem \ref{b787}, we should notice from the van Kampen theorem that the fundamental groups $ \pi_1( \Sigma_{p/1} (K,K'))$ and $\pi_1( \Sigma_{1/q} (K,K'))$ are presented by
\begin{equation}\label{b6522}
\langle \ x_1, \dots, x_{2g}, \mathfrak{m}
\ x_1', \dots, x_{2g'}' \mid r_1,r_2, \dots, r_{2g}, r_1', \dots, r_{2g'}', \ \ \ r_{2g+1}:= \mathfrak{l} \mathfrak{m}^p (\mathfrak{l}')^{-1}\ \ \rangle, \end{equation}
\begin{equation}\label{b652200}
\langle \ x_1, \dots, x_{2g}, \mathfrak{m}
\ x_1', \dots, x_{2g'}' , \mathfrak{m}'\mid r_1,r_2, \dots, r_{2g}, r_1', \dots, r_{2g'}', \ \ \ r_{\dagger}:= \mathfrak{m} \mathfrak{l}^q (\mathfrak{l}')^{-1}, r_{\star}:= \mathfrak{m}' \mathfrak{l}^{-1}\ \ \rangle. \notag \end{equation}
Here, in \eqref{b6522}, we identify $ \mathfrak{m} $ with $ \mathfrak{m} '$. Define $w_{i}$ to be $\prod_{j=1}^i [ u_{2j-1}^{\flat}, u_{2j}^{\flat} ]$, and
\[ W_{i} := w_{i-1} \rho_{2i-1}w_{i-1}^{-1} \cdot
(w_{i-1} u_{2i-1}^{\flat}) \rho_{2i }(w_{i-1} u_{2i-1}^{\flat})^{-1}
\cdot (w_{i} u_{2i}^{\flat}) \rho_{2i-1}^{-1}(w_{i} u_{2i}^{\flat})^{-1} \cdot w_{i} \rho_{2i }^{-1}w_{i}^{-1} . \]
Likewise, we also define words $w_i'$ and $ W_i'$. We consider the two words,
\[ W_{p/1}^{K,K'} := (\Pi_{i=1}^g W_i) \cdot \rho_{2g+1} \cdot (\Pi_{i=1}^{g'} W_i')^{-1} (\mathfrak{m} \rho_{2g+1 }^{-1} \mathfrak{m}^{-1 } ) , \]
\[ W_{1/q}^{K,K'} := (\Pi_{i=1}^g W_i )\cdot \rho_{\star}^{-1} \cdot (\mathfrak{m}' \rho_{\dagger}(\mathfrak{m}')^{-1} ) \cdot (\Pi_{i=1}^{g'} W_i' )\cdot (\mathfrak{l}' \rho_{\star} (\mathfrak{l}')^{-1}) \cdot \rho_{\dagger}^{-1} .\]
\begin{thm}\label{b787}
Then, $ W_{p/1}^{K,K'} $ and $W_{1/q}^{K,K'} $ are taut identities with respect to the presentations \eqref{b6522} of $\pi_1( \Sigma_{p/1} (K,K')) $ and $\pi_1( \Sigma_{1/q} (K,K'))$, respectively.
\end{thm}
\begin{proof}
An immediate computation gives $ \psi(W_i)= w_{i-1} \mathfrak{m} [ u_{2i-1}^{\sharp}, u_{2i}^{\sharp} ] \mathfrak{m}^{-1} w_i^{-1} $, so that $ \psi(\Pi_{i=1}^g W_i)=\mathfrak{m} \Pi_{i=1}^g [u_{2i-1}^{\sharp}, u_{2i}^{\sharp} ] \mathfrak{m}^{-1} w_g^{-1} $. Then, $ W_{p/1}^{K,K'} $ and $W_{1/q}^{K,K'} $ turn out to be identities by \eqref{b6522}. Furthermore, by the definition of $W^{K,K'}_{\bullet}$, we verify that $W^{K,K'}_{\bullet}$ are taut.
\end{proof}
As a corollary, if $K'$ is the unknot, we have the complex $C_*(\widetilde{M};\Z) $, where $M$ is the 3-manifold, $M_{p/1}(K)$, obtained by $p/1$-surgery of $S^3$ along $K$:
\begin{cor}\label{co14}
If $M:=M_{p/1}(K) $ satisfies Assumption ($\dagger$), then the boundary maps $\partial_2$ and $\partial_3 $ in the associated complex $C_*(\widetilde{M};\Z) $ in \eqref{b78} are given by the following matrix presentations:
\[ \left(
\begin{array}{cc}
\Bigl\{ \mathfrak{m} \frac{\partial u_{i}^{\sharp}}{\partial x_{j}} - \frac{\partial u_{i}^{\flat}}{\partial x_{j}} \Bigr\}_{1 \leq i,j \leq 2g} & \Bigl\{ 1- \frac{\partial \mathfrak{l}}{\partial x_{j}}  \Bigr\}^{\rm transpose}_{1 \leq j \leq 2g}\\
\Bigl\{ \frac{\partial \mathfrak{l} }{\partial x_{j}} \mathfrak{m}^p \Bigr\}_{1 \leq j \leq 2g}
& \mathfrak{l} \frac{\partial \mathfrak{m}^p }{\partial \mathfrak{m}} \\
\end{array}
\right) , \]
\[ \partial_3 (s)= (1- \mathfrak{m})b_{2g+1}+\sum_{i=1 }^g ( w_{i-1} -w_i u_{2i}^{\flat})b_{2i-1} + ( w_{i-1} u_{2i-1}^{\flat} - w_{i})b_{2i}. \]
\end{cor}
\begin{rem}\label{rem14} We give a comparison to Theorem 3.9 in \cite{MP}. The authors give an expression of the chain complex $C_*(\widetilde{M};\Z) $, where $M=M_{0/1}(K) $. However, the numbers of basis of $C_3,C_2,C_1$ are $2, c+1,c$, respectively, where $c$ is the crossing number of $K$, while
those in Corollary \ref{co14} are fewer.
\end{rem}

Let us recall the cabling conjecture, which predicts that if $K$ is not a cabling knot, then $ M_0(K) $ is irreducible; this conjecture has been proven for some classes of knots. Since $ \pi_1( M_0(K) )$ is of infinite order, it is fair to say that most $ M_0(K) $ satisfy Assumption ($\dagger$). Incidentally, it is a problem for the future to clarify a taut identity for the $(p/q)$-surgery for any $p/q \in \Q$.

\subsection{Branched covering spaces of $S^3$ branched over a knot}
\label{defS32}
Take a knot $K$ in $S^3$, and $ d \in \N$. In this subsection, we will give a taut identity of $ \pi_1(B_K^d)$, where we mean by $B_K^d$ the $d$-fold cyclic covering space of $S^3$ branched over $K$. We should remark the fact that, if $K$ is a prime knot and $\pi_1 ( B_K^d)$ is of infinite order, then $B_K^d $ is aspherical and therefore admits Assumption ($\dagger$). Let $p: E_K^d \ra S^3 \setminus K$ be the $d$-fold cyclic covering. For $ k \in \Z/d $, let $ x_i^{(k)}$ be a copy of $ x_i$ and $u_{i,k}^{\sharp}$ be the word obtained by replacing $x_i$ with $ x_i^{(k)}$ in the word $u_i^{\sharp}$. We similarly define the word $ u_{i,k}^{\flat }$. Then, by using the Reidemeister-Schreier method (see, e.g., \cite[Proposition 3.1]{Kab}), it follows from presentation \eqref{oo4656} that $ \pi_1(E_K^d )$ is presented by
\begin{equation}\label{b17} \langle \ x_1^{(k)}, \dots, x_{2g}^{(k)}, \overline{\mathfrak{m}} \ \ \ ( k \in \Z/d) \mid \overline{\mathfrak{m}} u_{i,k}^{\sharp} \overline{\mathfrak{m}}^{-1} \ ( u_{i,k+1}^{\flat })^{-1} \ \ \ (1 \leq i \leq 2g, \ k \in \Z/d) \ \ \rangle. \end{equation}
Since $B_K^d $ is obtained from $E_K^d$ by attaching a solid torus which annihilates the meridian $\overline{\mathfrak{m}} $, $ \pi_1(B_K^d)$ is presented by the quotient of $ \pi_1(E_K^d)$ subject to $\overline{\mathfrak{m}} =1$; that is,
\begin{equation}\label{b7} \pi_1(B_K^d) \cong \langle \ x_1^{(k)}, \dots, x_{2g}^{(k)} \ \ \ ( k \in \Z/d) \mid r_{i,k} := u_{i,k}^{\sharp} ( u_{i,k+1}^{\flat })^{-1}\ \ \ (1 \leq i \leq 2g, \ k \in \Z/d) \ \ \rangle. \end{equation}
Let $F $ be the free group $ \langle \ x_1^{(k)}, \dots, x_{2g}^{(k)} \ \ \ ( k \in \Z/d) \mid \rangle $. From \eqref{b652}, we should notice that $[ u_{1,k}^{\flat}, u_{2,k}^{\flat} ] \cdots [ u_{2g-1,k}^{\flat}, u_{2g,k}^{\flat} ] = [ u_{1,k}^{\sharp}, u_{2,k}^{\sharp} ] \cdots [ u_{2g-1,k}^{\sharp}, u_{2g,k}^{\sharp} ] \in F$ for any $ k \in \Z/d$.

Similarly to \S \ref{defS33}, we will give an identity with respect to the presentation \eqref{b7}. For $1 \leq i \leq g, 1\leq k \leq d $, define $w_{i,k}= \prod_{j=1}^i [ u_{2j-1,k+1}^{\flat}, u_{2j,k+1}^{\flat} ]$, and
\[ W_{i,k} = w_{i-1,k} \rho_{2i-1,k}w_{i-1,k}^{-1} \cdot
(w_{i-1,k} u_{2i-1,k+1}^{\flat}) \rho_{2i,k }(w_{i-1,k} u_{2i-1,k +1}^{\flat})^{-1}
\cdot \]
\[ \ \ \ \ \ \ \ \ \ \ \ \ \ (w_{i,k} u_{2i,k+1}^{\flat}) \rho_{2i-1,k }^{-1}(w_{i,k} u_{2i,k+1}^{\flat})^{-1} \cdot w_{i,k} \rho_{2i,k }^{-1}w_{i,k}^{-1} . \]
\begin{prop}\label{b7802}
Define $W$ to be $\Pi_{k=1}^d W_{1,k} W_{2,k }\cdots W_{g,k} $, by the above equality in $F$. Then, $W$ is a taut identity. In particular, if $B_K^d$ satisfies Assumption ($\dagger$), the associated complex in \eqref{b78} is isomorphic to the cellular chain complex of the universal cover of $B_K^d $.
\end{prop}
\begin{proof}
Direct calculation gives $ \psi( W_{i,k})= w_{i-1,k} [u_{2i-1,k}^{\sharp}, u_{2i,k}^{\sharp} ] w_{i,k}^{-1}$, which deduces
$$ \psi ( \Pi_{i=1}^g W_{i,k})= ( \Pi_{i=1}^g[ u_{2i-1,k}^{\sharp}, u_{2i,k}^{\sharp} ] )w_{g,k}^{-1} =
\Pi_{i=1}^g[ u_{2i-1,k}^{\sharp}, u_{2i,k}^{\sharp} ] ( \Pi_{i=1}^g[ u_{2i-1,k+1}^{\flat}, u_{2i,k+1}^{\flat} ])^{-1} . $$
Thus, $W$ turns out to be an identity. Furthermore, since we can verify that $W$ is taut by the definition of $W$, Remark \ref{bb2} readily leads to the latter part.
\end{proof}
\begin{exa}\label{exa78}
Let $K$ be the figure-eight knot. It can be verified that the presentation \eqref{oo4656} can be written as
$$ \langle \ x_1, x_2 , \mathfrak{m} \mid \mathfrak{m} x_1x_2 \mathfrak{m}^{-1} = x_1 , \ \mathfrak{m} x_2 x_1x_2 \mathfrak{m}^{-1} = x_2 \ \ \rangle. $$
Thus, by \eqref{b7}, we have
$$ \pi_1(B_K^d) \cong \langle \ x_1^{(i)}, x_2^{(i)} \ \ (1 \leq i\leq d) \mid x_1^{(i)}x_2^{(i)} = x_1^{(i+1)} , \ \ x_2^{(i)} x_1^{(i)}x_2^{(i)} = x_2^{(i+1)}\ \ (1 \leq i\leq d) \ \ \rangle. $$
Annihilating $x_2^{(i)} $ by using the relation $x_1^{(i)}x_2^{(i)} = x_1^{(i+1)}$, we have
$$ \pi_1(B_K^d) \cong \langle \ x_1^{(1)}, \dots, x_1^{(d)} \mid (x_1^{(i)})^{-1}(x_1^{(i+1)})^2
(x_1^{(i+2)})^{-1}x_1^{(i+1)} \ \ (1 \leq i\leq d) \ \rangle. $$
This isomorphism coincides exactly with the result in \cite[Page 963]{KKV}.

Likewise, we can verify that some groups, called ``cyclically presented groups" in \cite{KKV} and references therein, are isomorphic to $\pi_1(B_K^d) $ for some $K$ and $d$.
\end{exa}
\begin{rem}\label{rea78} As the referee points out, it is reasonable to hope that Proposition \ref{b7802} is true without Assumption ($\dagger$).
In fact, as seen in \cite{Sie1}, given a Heegaard diagram, we can construct a ``squashing map" and a taut identity compatible with the complex \eqref{b78}. Thus, it is a conjecture that we can find an appropriate Heegaard diagram of $ B_K^d$ such that the associated taut identity is equal to the above $W$. 
\end{rem}

\subsection{0-Surgery-like spaces from branched covering spaces of $S^3$}
\label{defS34}
Using the notation in the preceding subsection, we can examine the 3-manifold obtained by the 0-surgery on the knot $p^{-1}(K) \subset B^d_K$. 
The 0-surgery appears in the topic of the concordance group including the Casson-Gordon invariant \cite{CG}. More precisely, regarding the boundary of $ E_K^d$ as a knot in $ B_K^d$, we consider the 3-manifold obtained by 0-surgery on the knot in $ B_K^d$. Notice from \eqref{b17} that the fundamental group canonically has a group presentation
\begin{equation}\label{b7768} \langle \ x_1^{(k)}, \dots, x_{2g}^{(k)} \ \ \ ( k \in \Z/d) , \ \overline{\mathfrak{m}} \mid r_i^{(k )} \ \ \ ( i \leq 2g, \ k \in \Z/d), \ \ \Pi_{i=1}^g[ u_{2i-1,1}^{\flat}, u_{2i ,1}^{\flat} ] \rangle. \end{equation}
Let $\mathfrak{ l }^{(k)}:= \Pi_{i=1}^g[ u_{2i-1,k}^{\flat}, u_{2i ,k}^{\flat} ]$, and consider an analogous presentation
\begin{equation}\label{b776} \langle \ x_1^{(k)}, \dots, x_{2g}^{(k)} \ \ \ ( k \in \Z/d) , \ \overline{\mathfrak{m}} \mid r_i^{(k )} \ \ \ ( i \leq 2g, \ k \in \Z/d), \ \ r_{\ell}:=
\mathfrak{ l }^{(1)} \mathfrak{ l }^{(2)} \cdots \mathfrak{ l }^{(d)} \rangle. \end{equation}
Similarly to \S \ref{defS32}, we can construct an identity. For $ i \leq 2g, k \leq d $, define
$z_{k}= \mathfrak{ l }^{(1)} \mathfrak{ l }^{(2)} \cdots \mathfrak{ l }^{(k)} $ and
\[ W_{i,k} = z_k w_{i-1,k} \rho_{2i-1,k}w_{i-1,k}^{-1} z_k^{-1} \cdot
( z_k w_{i-1,k} u_{2i-1,k+1}^{\flat}) \rho_{2i,k }( z_k w_{i-1,k} u_{2i-1,k +1}^{\flat})^{-1}
\cdot \]
\[ \ \ \ \ \ \ \ \ \ \ \ \ \ ( z_k w_{i,k} u_{2i,k+1}^{\flat}) \rho_{2i-1,k }^{-1}( z_k w_{i,k} u_{2i,k+1}^{\flat})^{-1} \cdot z_k w_{i,k} \rho_{2i,k }^{-1}w_{i,k}^{-1} z_k^{-1} . \]
In the usual way, we can easily show the following:
\begin{prop}\label{b7812}
Define $W$ to be $(\Pi_{k=1}^d \Pi_{i=1}^g W_{i,k} )\cdot \rho_{\ell} \cdot (\mathfrak{m} \rho_{\ell}^{-1} \mathfrak{m}^{-1} ) $. Then, $W$ is a taut identity. In particular, Remark \ref{bb2} ensures
that if the fundamental group of a closed 3-manifold satisfying Assumption $(\dagger)$ is isomorphic to
\eqref{b776}, then
the cellular chain complex of the universal cover is isomorphic to the complex \eqref{b78}.
\end{prop}

\subsection{Some Seifert fibered spaces over $S^2$}
\label{defS35}
In the last subsection, we will discuss some of the Seifert fibered spaces and Brieskorn manifolds. The theorem of Scott \cite{Sc} shows that the homeomorphism types of such spaces with infinite $\pi_1$ can be detected by the fundamental groups; thus, the spaces satisfy Assumption ($\dagger$).

Let us state Proposition \ref{b7813216}. Take integers $a_1, \dots, a_{n+1}$ with $ a_i \geq 2$, and $\epsilon_1 ,\dots, \epsilon_n \in \{ \pm 1\}$. Let $M$ be a Seifert fibered space of the form
$$ \Sigma (0; (0,1), (a_1, \epsilon_1 ),(a_2, \epsilon_2),\dots, (a_{n}, \epsilon_n) , (a_{n+1}, 1) ).$$
Then, as is classically known, the fundamental group has the presentation
$$ \langle \ x_1, \dots, x_{n+1}, h \mid h x_i h^{-1} x_{i}^{-1}, \ \ x_{i}^{a_i}h^{\epsilon_i} \ \ (i \leq n) \ \ \ \ , x_{n+1}^{a_{n+1}}h ,\ \ x_1 \cdots x_{n+1}\ \ \rangle. $$
Furthermore, let us consider a group $G$ with the presentation
\begin{equation}\label{b7796} \langle \ x_1, \dots, x_{n} \mid r_i:= (x_i x_{i+1} \cdots x_n x_{ 1}\cdots x_{i-1} )^{-a_{n+1}} x_{i}^{\epsilon_i a_i} \ \ \ \ ( i \leq n) \ \ \rangle. \end{equation}
We can easily check that the correspondence $ x_i \mapsto x_i, \ \ x_{n+1} \mapsto (x_1 \cdots x_n)^{-1}, \ \ h\mapsto x_1^{\epsilon_1 a_1 } $ gives rise to a group isomorphism $ \pi_1(M) \cong G.$ Therefore, we shall define a taut identity on the presentation \eqref{b7796}:
\begin{prop}\label{b7813216}
Suppose that $\pi_1(M)$ is of infinite order. Define $W$ to be
$$ \rho_1( x_1^{-1} \rho_1^{-1} x_1) \rho_2( x_2^{-1} \rho_2^{-1} x_2)\cdots
\rho_n( x_n^{-1} \rho_n^{-1} x_n) . $$
Then, $W$ is a taut identity of the presentation \eqref{b7796}.
\end{prop}
The proof is similar to the ones above, so we will omit the details.
\begin{rem}\label{bb132}
The taut identity when $n=2$ is presented in \cite[p. 127]{Sie1}. The paper does not mention the homeomorphism type of the associated 3-manifold; however, Proposition \ref{b7813216} implies that the homeomorphism type can be detected by a Seifert structure.
\end{rem}

Finally, let us turn to the topic of Brieskorn 3-manifolds. Choose integers $a,b,p, q,m \in \Z$ and $\varepsilon \in \{ \pm 1 \} $ satisfying $ap+ bq =1$ and $p,q,m >1.$ We will focus on the Brieskorn 3-manifold of the form,
$$M:=\Sigma(p,q,m pq + \varepsilon ):= \{ (x,y,z) \in \mathbb{C}^3 \mid x^p+y^q+z^{m pq + \varepsilon } =0, \ \ |x|^2+|y|^2+|z|^2= 1\ \ \}, $$
which is an Eilenberg-MacLane space if $ 1/p+1/q+1/( m pq + \varepsilon ) <1$. The manifold is known to be homeomorphic to a $3$-manifold obtained from $(\varepsilon /m)$-surgery on the $(p,q)$-torus knot $T_{p,q}$. Recall the presentation of $ \pi_1( S^3 \setminus T_{p,q} )$ as $ \pi_1(S^3 \setminus T_{p,q}) \cong \langle \ x,y \mid x^q =y^p \ \rangle $, and that the meridian $\mathfrak{m}$ and the preferred longitude $\mathfrak{l}$ are identified with $x^a y^b$ and $(x^a y^b )^{-pq } x^q$, respectively. Therefore, $\pi_1(M)$ admits a genus-two Heegaard decomposition and has the group presentation,
\begin{equation}\label{kk4}\pi_1(M) \cong \langle \ x,y \mid r_1 :=x^{qm} (x^a y^b)^{-mpq - \varepsilon} , \
r_2 := (x^a y^b)^{mpq + \varepsilon} y^{-p} x^{-qm-q} \ \ \rangle. \end{equation}
Likewise, we can show the following result:
\begin{prop}\label{b78132} Suppose $ 1/p+1/q+1/( m pq + \varepsilon ) <1$ as above.
Then the following word is a taut identity of the presentation \eqref{kk4}.
\[ \rho_1\rho_2^{-1} \rho_1^{-1} (x^{qm}y^{-p} x^{-qm-q} \rho_2 x^{qm+q}y^{p} x^{-qm}) . \]
\end{prop}

\section{First application to the linking forms of branched covers}
\label{defS41}
\subsection{Review of the linking form and a theorem}
\label{defS412}
Here, we will review the linking form of $M$ for a closed 3-manifold $M$ with $H_*(M;\Q) \cong H_*(S^3 ;\Q)$. Considering the short exact sequence
\begin{equation}\label{b1} 0 \lra \Z \lra \Q \lra \Q/\Z \lra 0, \end{equation}
we can easily check that the Bockstein maps
$$ \beta: H_2( M;\Q/\Z) \cong H_1 (M ;\Z), \ \ \ \ \ \beta: H^1( M;\Q/\Z) \cong H^2 (M ;\Z), $$
are isomorphisms from the long exact homology sequences. Let $\PD_M^{\Z}$ be the Poincar\'{e} duality on the integral (co)-homology. We denote by $\Omega$ the composite map defined by setting
$$ H_1(M;\Z) \xrightarrow{\ \PD_M^{\Z} \ } H^2(M;\Z) \xrightarrow{ \ \beta^{-1} \ } H^1 (M;\Q/ \Z) \xrightarrow{ \ \rm ev \ } \Hom( H_1(M;\Z);\Q/\Z) ,$$
where the last map is the Kronecker evaluation map. Then, {\it the linking form of $M$} is
$$ \lambda_{M}: H_1(M;\Z) \times H_1(M;\Z) \lra \Q/\Z $$
defined by $ \lambda_M (a, b)= \Omega (a )(b)$. This bilinear map is known to be symmetric and non-singular. This definition goes back to Seifert \cite{Sei}, and the form has sometimes appeared in the study of algebraic surgery theory (see, e.g., \cite{Wall}) and the concordance groups of knots \cite{CG}. Recently, the linking form of $M$ can be computed in terms of Heegaard splittings \cite{CFH}.

Of particular interest to us is an application to the Casson-Gordon invariant \cite{CG} and a procedure for computing $\lambda_{M} $ in another way. In what follows, let $B_K^d $ be the $d$-fold cyclic covering space of $S^3$ branched over a knot $K$. In the context of the invariant, the linking form of $B_K^d $ plays an important role: more precisely, it is important to calculate metabolizers of the form; see, e.g., \cite{CG}.

Now let us give a matrix presentation of the homology $H_1(B_K^d;\Z)$ and state the main theorem. Choose a Seifert surface $\Sigma $ of $K$ whose genus is $g$, as in \S \ref{defS33}. Then, we have the Seifert form $\alpha : H_1(\Sigma;\Z) \otimes H_1(\Sigma;\Z) \ra \Z$; see \cite[Chapter 6]{Lic} for the definition. Let $J$ be the inverse matrix $ (V - \ \! \! ^t V )^{-1}$, where $\mathrm{det}(V - \ \! \! \! ^t V )=1$ is known (see \cite[Theorem 6.10]{Lic}). The matrix presentation is often written as $V \in \mathrm{Mat}(2g\times 2g;\Z)$ and is called {\it the Seifert matrix}. Consider the following matrices of size $(2g d \times 2gd)$:
$$ A:= \left(
\begin{array}{ccccc}
-V & 0 & \cdots &0& ^t V \\
^t V & -V & \cdots &0 & 0 \\
0 & ^t V &\ddots & \vdots & \vdots \\
\vdots & \vdots &\ddots & \ddots & \vdots \\
0 & 0 &\cdots & ^t V & -V
\end{array}
\right), \ \ \ \ \ \ \ \ B:= \left(
\begin{array}{ccccc}
0 & 0 & \cdots &0& J^t V \\
J^tV & 0 & \cdots &0 & 0 \\
0 & J^tV &\ddots & \vdots & \vdots \\
\vdots & \vdots &\ddots & \ddots & \vdots \\
0 & 0 &\cdots & J^tV & 0
\end{array}
\right) , $$
which appear in \cite[Page 494]{Tro}.
As is known (see \cite{Sei,Tro} or \cite[Theorem 9.7]{Lic}), the first homology $H_1(B_K^d;\Z)$ is isomorphic to the cokernel of $A$, i.e., $H_1(B_K^d;\Z) \cong \Z^{2gd}/ \ \! ^t \! A\Z^{2gd}. $ In particular, $\mathrm{det}( A)\neq 0$ if and only if $H_1 (B_K^d ;\Q) \cong 0 $. The linking formula of $ B_K^d $ can be algebraically formulated in the above notation as follows:
\begin{thm}\label{bbbn1}
Suppose that $B_K^d $ satisfies Assumption ($\dagger$) and $H_1 (B_K^d ;\Q) \cong 0 $. Then, the matrix multiplication $B : \Z^{2gd} \ra \Z^{2gd} $ induces an isomorphism
$ \mathcal{B}: \Z^{2gd}/ \ \! ^t \! A\Z^{2gd} \ra \Z^{2gd}/ \ \! ^t \! A\Z^{2gd}$ and the linking form $ \lambda_{B_K^d }$ of $B_K^d $ is equal to the form,
\begin{equation}\label{b212} \Z^{2gd}/ \ \! ^t \! A\Z^{2gd} \times \Z^{2gd}/ \ \! ^t \!A\Z^{2gd} \lra \Q/\Z; \ \ \ \ \ \ (v,w) \longmapsto
\ ^t v \mathrm{adj}(A) \!\! \ ^t \mathcal{B}^{-1} w/\Delta . \end{equation}
Here, $\mathrm{adj}(A) $ is the adjugate matrix of $A$, and $\Delta $ is the order $| H_1 (B_K^d ;\Z)|\in \mathbb{N}. $
\end{thm}
This statement is implicitly connoted in \cite[Satz I]{Sei} and \cite[p. 496]{Tro} \footnote{To be precise, the original statements implicitly claim that the linking form $ \lambda_{B_K^d }$ is equal to the matrix presentation $ B \mathrm{adj}(A) $ up to isomorphisms. However, for applications to the Casson-Gordon invariants, we should describe the linking form from a basis of $ H_1(B_K^d;\Z)$.}
; however, there is no complete proof for this statement in the literature.

Here, let us make a few remarks. Whereas the matrix $\mathrm{adj}(A) \!\! \ ^t \mathcal{B}^{-1} $ is not always symmetric, the quotient on $\Z^{2gd}/ \ \! ^t \! A\Z^{2gd} $ is symmetric. Next, the second condition of $ H_1 (B_K^d ;\Q) \cong 0 $ is not so strong: indeed, according to \cite[Corollary 9.8]{Lic}, if any $d$-th root of unity is not a zero point of the Alexander polynomial of $K$ (e.g., the case $d$ is a prime power), then $H_1 (B_K^d ;\Q) \cong 0 $. Furthermore, as the proof and Remark \ref{rea78} imply, one may hope that the theorem is true even if we drop the condition of Assumption ($\dagger$).

\begin{proof}[Proof of Theorem \ref{bbbn1}]
It is known \cite[Lemma 2.5]{CFH} that the linking form can be formulated in the terminology of cohomology as
\begin{equation}\label{b22} \lambda_{M} (a,b)= \langle (\beta^{-1} \circ \PD_M^{\Z} ) (a) \smile \PD_M^{\Z} (b), [M]\rangle . \end{equation}
Here, $ \smile$ is the cup-product $H^1(M;\Q/\Z ) \otimes H^2(M;\Z ) \ra H^3(M;\Q/\Z )$.

Let $M$ be $ B_K^d$, and let $R$ be one of $\Z, \Q$ or $\Q/\Z$ as trivial coefficients. Let $\varepsilon: \Z[\pi_1(M )]\ra \Z$ be the augmentation map. Then, as is known (see \cite[Proposition 4.1]{Tro}), by choosing a Seifert surface, the integral matrices $\{ \varepsilon (\frac{\partial u_i ^{\sharp}}{\partial x_j} )\} _{1 \leq i,j \leq 2g} $ and $\{ \varepsilon (\frac{\partial u_i ^{\flat}}{\partial x_j} )\} _{1 \leq i,j \leq 2g} $ are equal to $V$ and $^t V $, respectively. Let us identify the complex $C^* (M;R) $ in the coefficients $R$ with $ C^* (\widetilde{M} ;\Z) \otimes_{\Z[\pi_1(M )]} R$ via $\varepsilon.$ Then, by presentation \eqref{b7}, the complex $C^* (M;R) $ reduces to
\begin{equation}\label{b2442}0 \ra C^0(M;\Z)\stackrel{0}{\lra} C^1 (M;R) \xrightarrow{ \ A \ } C^2 (M;R) \stackrel{0}{\lra} C^3 (M;R) \ra 0. \end{equation}
If $R=\Q$, the matrix $A$ is an isomorphic because of $H^* (M;\Q) \cong H^* (S^3 ;\Q)$. Therefore, from the definition of the Bockstein inverse map $ \beta^{-1}: C^2 (M;\Z) \ra C^1 (M;\Q/\Z) $ is identified with $ \Z^{2gd} \ra (\Q/\Z)^{2gd}; v \mapsto \mathrm{adj}(A)v/ \Delta $.

Meanwhile, from the formula for the identity $W$ in Proposition \ref{b7802} and the formula \eqref{b8w}, the cup-product $ \smile: C^1 (M;\Q/\Z) \times C^2 (M;\Z) \ra C^3 (M;\Q/\Z) \cong \Q/\Z $ is considered to be $ (\Q/\Z)^{2gd} \times \Z^{2gd} \ra \Q /\Z; (v,w) \mapsto ^t \!\! v B w $. The Poincar\'{e} duality ensures the non-degeneracy of the cup product on cohomology.
In particular, the desired induced map $\mathcal{B} $ is an isomorphism, and
is identified with the duality $ H_1( M;\Z) \cong H^2( M;\Z) $, where $H^2( M;\Z) $ is canonically regarded as 
$\mathrm{Coker}(A)= \Z^{2gd}/A\Z^{2gd} $ by \eqref{b2442}. Hence, upon the identification $ H^2( M;\Z) \cong H_1( M;\Z) \cong \Z^{2gd}/^tA\Z^{2gd}$, the formula \eqref{b22} immediately implies that the linking form is equal to the required \eqref{b212}.
\end{proof}

\subsection{Example computations }\label{defS12}
It is easier to quantitatively compute kernels rather than cokernels. Let us examine Corollary \ref{bbbn2} below. Let $ \Ker (A)_{\Z/\Delta} $ be $ \{v \in (\Z/\Delta \Z)^{2gd} \mid Av=0 \in (\Z/\Delta \Z)^{2gd}\ \} $. Consider the linear map
$$ \Z^{2gd}/A\Z^{2gd} \lra \Ker (A)_{\Z/\Delta}; \ \ v \longmapsto \mathrm{adj}(A) v .$$
This map is an isomorphism if $|\Delta | \neq 0 $: in fact, with a choice of the section $\mathfrak{s}: \Z^{2gd}/A\Z^{2gd} \ra \Z^{2gd}$, the inverse map is defined by $w \mapsto (A \mathfrak{s}(w))/\Delta$. In summary, from Theorem \ref{bbbn1}, we immediately have the following:
\begin{cor}\label{bbbn2}
Let $\Delta$, $A,B$ and $\mathrm{adj}(A) $ be as in Theorem \ref{bbbn1}. Under the supposition in Theorem \ref{bbbn1}, the linking form $ \lambda_{B_K^d }$ of $B_K^d $ is isomorphic to the bilinear form
$$\Ker (A)_{\Z/\Delta} \times \Ker (A)_{\Z/\Delta} \lra \Q/\Z; \ \ \ \ \ \ (v,w) \longmapsto
\ ^t \mathfrak{s}(v) ^t A B w/\Delta^2 . $$
\end{cor}

\begin{exa}\label{yy1}
Let $p,q,r \in \Z$ be odd numbers. Let $K$ be the Pretzel knot $P(p,q,r)$. When $d=2$, the branched cover $B_K^2$ is known to be a Seifert fibered space of type $\Sigma(p,q,r)$ over $S^2.$ Furthermore, we can choose a Seifert matrix of the form $V= \frac{1}{2}\left( \begin{array}{cc} p+q & q +1 \\ q-1 & q+r \end{array} \right)$, and $\Delta= pq + q r + rp$; see \cite[Example 6.9]{Lic}.

First, consider the case where $p,q,r$ are relatively prime. Then, $\Ker(A)$ is generated by $(-r - q, q, -r - q, q)$, and we can easily verify that the linking form equal to $2 (q+r)/ \Delta $.

However, if $p,q,r$ are not relatively prime, $\Ker(A)$ and the linking form are complicated. For example, if $(p,q,r)=(p,-p,p)$, then $\Ker(A) \cong (\Z/p)^2$ is generated by $(0,p,0,p)$ and $(p,0,p,0)$; the linking matrix is equal to $\frac{2}{p}\left( \begin{array}{cc} 0& 1 \\ 1 & 0 \end{array} \right) $. Meanwhile, if $(p,q,r)=(p,p,p)$ and $ p$ is not divisible by 3, then $\Ker(A) \cong \Z/p \oplus \Z/3p$ possesses a basis, $v= (0,3p,0,3p)$, $w=(3p +p^2, p^2,3p+p^2,p^2)$. 
Hence, $\left( \begin{array}{cc} \mathrm{lk}(v,v)& \mathrm{lk}(v,w) \\ \mathrm{lk}(w,v) & \mathrm{lk}(w,w) \end{array} \right)$ can be computed as $\frac{2}{p}\left( \begin{array}{cc} 2& 1 \\ 1 & 2 \end{array} \right) $.

In a similar way, we can compute many linking forms of $d$-fold branched covering spaces for small $d$ with the help of a computer program.
\end{exa}

\section{Second application to Dijkgraaf-Witten invariants}
\label{defS42}
As another application, we develop procedures of computing some Dijkgraaf-Witten invariants in terms of identities.

We start by reviewing the Dijkgraaf-Witten invariant \cite{DW}.
Let $G$ be a finite group, $A$ a commutative ring, and $ \psi$ a
group 3-cocycle of $G$.
Denoting by $B G$ an Eilenberg-MacLane space of type $(G, 1)$, we have
a classifying map $\iota: M \hookrightarrow B\pi_1(M)$ uniquely up to homotopy.
Then, as is known, $\psi$ can be regarded as a 3-cocycle of $H^3(BG;A) $, and any group homomorphism $f : M \ra G$ canonically gives rise to the composite
$$ \iota^* \circ f^*:H^*(B G;A ) \ra H^*( B\pi_1(M);A )= H^*( \pi_1(M);A ) \ra H^*( M ;A ) .$$
Then, the {\it Dijkgraaf-Witten invariant of $M$} is defined as a formal sum in the group ring $\Z[A]$ by setting
$$ \DW_{\psi}(M) := 
\sum_{f \in \mathrm{Hom}(\pi_1(M),G)}\langle \iota^* \circ f^* (\psi), [M] \rangle \in \Z [A] . $$
Although the definition seems rather simple or direct, it is not easy to compute $ \DW_{\psi}(M) $ except in the case where $G$ is abelian, because it is not trivial to explicitly express $[M]$ and $f^*$ (however, see \cite{DW,W2} for the abelian case and \cite{No2} for a partially non-abelian case).
To the knowledge of the author, there are few examples of Dijkgraaf-Witten invariants when $G$ is non-abelian.

This section develops a method for computing the invariants, and gives non-abelian examples.
First, for simplicity, we now restrict on the case $\psi = \gamma \smile \delta$ for some
$\gamma \in H^1(G ;A)$ and $\delta \in H^2(G ;A)$. Take a group homomorphism $f:\pi_1(M) \ra G$ and a group presentation
$G= \langle y_1, \dots, y_n \ \mid \ s_1, \dots, s_\ell \rangle $. Then, as in \eqref{b78}, we have a commutative diagram:
$${\normalsize
\xymatrix{ C_*(M;A): \ar[r]^{ \partial_3}& A \otimes \Z[\pi_1(M)]^m \ar[r]^{ \partial_2}\ar[d]_{f_*}& A \otimes \Z[\pi_1(M)]^m \ar[d]_{f_*} \ar[r]^{ 0} & A \otimes \Z[\pi_1(M)] \ar[d]_{f_*}\\
C_*(G;A): \ar[r]^{ \partial_3'} &A \otimes \Z[G]^\ell \ar[r]^{ \partial_2'} & A \otimes \Z[G]^n \ar[r]^{ \partial_1'=0} & A \otimes \Z[G].
}}$$
Here, the tensors are over $\Z[G]$, and $\partial_2'(b_i')=\sum_{k=1}^n [ \frac{\partial s_i}{\partial y_k}] a_k'$.
\begin{exa}\label{lll42}Suppose $p,q \in \mathbb{N}$ such that $(p,q)=1$. Let $A= G=\Z/p$, and $M $ be the lens space $L(p,q) $.
Then, as is known, $H^*(G;A)\cong \Z/p$, and we can choose appropriate generators $\alpha_i \in H^i(G;A)\cong \Z/p$ such that $\alpha_3 = \alpha_1 \smile \alpha_2$.
We fix a presentation $ G= \pi_1(M)= \langle x \ \mid \ | s:=x^p \rangle $. Then, the taut identity of $L(p,q)$ is known to be
$ W_{p,q}= s x^{-q} s^{-1} x^q $; see \cite{Sie1}.
Then, for $i \leq 3$, we can regard $\alpha_i$ as a map $\Z/p=C_*(M;\Z/p) \ra \Z/p$ that sends a generator to $1$. Then it follows from \eqref{b4} that the cup product
$\smile : H^1(L(p,q); \Z/p) \times H^2( L(p,q); \Z/p) \ra \Z/p$ is computed as $ (a,b) \mapsto q ab$. Moreover, for $a \in \Z/p $, if we define $f_a : \pi_1(M) \ra G$ by setting $x \mapsto a$, then $\Hom( \pi_1(M),G)$ is equal to $\{ f_a | a \in \Z/p \}$, and we can compute
$$\langle f_a^* ( \alpha_3), \iota_*[M]\rangle = \langle f_a^* ( \alpha_1 \smile \alpha_2), \iota_*[M]\rangle = \langle a \alpha_1 \smile a \alpha_2, \iota_*[M]\rangle=q a^2.$$
In conclusion,
$$ \DW_{\alpha_3} ( L(p,q)) = \sum_{ a \in \Z/p} 1\{ qa^2\} \in \Z[ \Z/p]. $$
In a similar way, if $M$ is another manifold such that the cohomology ring is known, we can compute $ \DW_{\alpha_3} ( M) $ for $G=\Z/p$. Comparing with \cite{DW,W2} as original computations, the above computation seems easier.
\end{exa}
\begin{exa}\label{lll6}
Let $m, n$ be natural numbers such that $m$ is relatively prime to $6n$. Let $G$ be
the non-abelian group of order $m^3$ which has a group presentation
\begin{equation}\label{PPP} \langle x,y, z \ | \ x^m,y^m ,z^m , s:= xzx^{-1}z^{-1},t:=yzy^{-1}z^{-1} , u:= zyxy^{-1}x^{-1}\rangle . \end{equation} 
The (co)-homology of $G$ is known (see, e.g., \cite{Lea}).
As a result, $H_1(G;\Z) \cong (\Z/m)^2$. 
Dually, the first cohomology $H^1(G;\Z/m) \cong (\Z/m)^2 $ is generated by the maps $\alpha$ and $\beta$ defined by $\alpha(x)=\beta(y)=1$ and $\alpha(y)=\beta(x)=0$.
Furthermore, the Massey product $\langle \alpha, \beta, \alpha \rangle $ 
and the product $\psi:= \beta \smile \langle \alpha, \beta, \alpha \rangle $ are known to be non-trivial. The equality $\psi =- \alpha \smile \langle \beta, \alpha,\beta \rangle $ is also known.
Since the cup product $C^1 \otimes C^1 \ra C^2$ is well described in \cite[\S 2.4]{Tro}, the Massey product $\langle \alpha, \beta, \alpha \rangle $ can be, by definition, regarded as the map $C_2(G;\Z/m) \ra \Z/m$ by setting
\begin{equation}\label{PPP1} x^m \mapsto 0, \ \ y^m \mapsto 0, \ \ z^m \mapsto 0, \ \ s \mapsto 0, \ \ t \mapsto 0,\ \ u\mapsto 2.\end{equation}

On the other hand, for simplicity, we specialize to 
the Seifert manifolds of type $M_{m,n} := \Sigma (0, (1,0),(m,1), (m, -1), (n,-1))$ over $S^2$,
whose fundamental groups are presented by
$$ \pi_1(M_{m}) = \langle x_1,x_2 \ | \ r_1:=x^m_1 (x_1^{-1}x_2^{-1})^n, \ r_2:= x_2^m(x_2^{-1}x_1^{-1})^n \rangle. $$
By Proposition \ref{b7813216}, the identity is $W:= r_2 x_2 r_2^{-1} x_2^{-1} r_1 x_1 r_1^{-1} x_1^{-1}$.
We further analyze the set $\Hom(\pi_1(M_{m,n}),G)$. For $a,b,c \in \Z/m$, consider the homomorphism $ f_{a,b,c}: \pi_1(M_{m,n}) \ra G$ defined by
$$ f_{a,b,c}(x_1):= x^a y^b z^c, \ \ \ f_{a,b,c}(x_2):= x^{-a} y^{-b} z^{-c+ab}.$$
It is not so hard to check the bijectivity of $(\Z/m)^3 \leftrightarrow \Hom(\pi_1(M_{m,n}),G)$ which sends $(a,b,c)$ to $f_{a,b,c}$. Then, the conclusion is as follows:
\begin{prop}\label{lll776}Let $ \psi$ be $\beta \smile \langle \alpha, \beta, \alpha \rangle \in H^3( G,\Z/m)$. Let $m \in \Z$ be relatively prime to $6n$. 
Then, upon the identification $ (\Z/m)^3 \leftrightarrow \Hom(\pi_1(M_{m,n}),G) $, the Dijkgraaf-Witten invariant is equal to
$$ \DW_{\psi}(M_{m,n})= \sum_{ (a,b,c) \in (\Z/m)^3}1\{ n\bigr(2abc- a(a-1)b(b-1) \bigl)\} \in \Z[\Z/m] .$$
\end{prop}
\begin{proof}
Recall from \eqref{b78} that the basis of $C_2(\widetilde{M_{m,n}}) \cong \Z[\pi_1(M_m)]^2 $ is denoted by $b_1,b_2$, where $b_i$ corresponds to the relator $r_i$.
We now analyse $ (f_{a,b,c})_* ( b_1)\in C_2(G;\Z/m)$. 
We can easily check that $f_{a,b,c}(r_i ) $ is transformed to $x^{am}y^{bm} z^{cm(m+1)/2 }$ by the above relators $s,t,u$.
Let us define $N_{b_i} \in \Z$ to be the numbers of applying $u$ when we transform $(f_{a,b,c})(r_i)$ by $x^{ma}y^{bm} z^{cm(m+1)/2 } $. 
Then, by \eqref{PPP1}, the pairing $\langle \langle \alpha, \beta, \alpha \rangle , \ (f_{a,b,c})_* ( b_1)\rangle$ is equal to $2N_{b_1}$.
From the definition of $N_{b_1} $, a little complicated computation can lead to
$$ N_{b_1}= \frac{m(m+1) ac}{2}+ (\sum_{i=1}^{m-1} \frac{ia(ia-1)}{2})+ nac - \frac{na(a-1)(b-1)}{2} \in \Z.$$
Since $m$ is relatively prime to $6n$, we can easily check the first and second terms to be zero modulo $m$.
Hence, using the above description of $W$ and the formula \eqref{b4}, we have
\[ \langle \psi,(f_{a,b,c})_*[M_{m,n}] \rangle =0+ b\cdot 2 N_{b_1} -0 \cdot N_{b_2}+0=b (2 nac - na(a-1)(b-1)) \in \Z/m,\]
which immediately leads to the conclusion.
\end{proof}
\end{exa}
\noindent
The above computation is relatively simple, since so are the presentations of $\pi_1(M)$ and $G$; however, a similar computation seems to be harder if $\pi_1(M)$ is complicated.

In contrast, we conclude this paper by suggesting another procedure of computing $\DW_{\psi}(M)$, which is
implicitly discussed in \cite[\S 4]{Nos}. Hereafter $\psi \in H^3(G;A)$ may be an arbitrary 3-cocycle.

Let $C_*^{\rm nh}(G;\Z) $ be the normalized homogenous complex of $G$, which is defined as the quotient $\Z$-free module of $\Z[G^{n+1}]$ subject to the relation $(g_0, \dots, g_n) \sim 0$ if $g_i = g_{i+1}$ for some $i$; see \cite[19 page]{Bro}.
Assume that we know an explicit expression of $\psi:G^4 \ra A$ as an element of $C^3_{\rm nh}(G,A) $.
When $* \leq 3$, we now define a chain map $c_*: C_n(\widetilde{M};\Z) \ra C_n^{\rm nh}(\pi_1(M);\Z)$ as follows. Let $c_0$ be the identity map. Let $A \in \Z[\pi_1(M)]$ be any element. Define $c_1(A x_i) := (A,Ax_i) $. If $r_i$ is expanded as $ x_{i_1}^{\epsilon_1} x_{i_2}^{\epsilon_2} \cdots x_{i_n}^{\epsilon_n} $ for some $ \epsilon_k \in \{ \pm 1 \}$, we define
$$c_2( A r_i) = \sum_{m: 1 \leq m \leq n}\epsilon_m ( A ,A x_{i_1}^{\epsilon_1} x_{i_2}^{\epsilon_2} \cdots
x_{i_{m-1}}^{\epsilon_{m-1}} x_{i_m}^{(\epsilon_m -1)/2} , A x_{i_1}^{\epsilon_1} x_{i_2}^{\epsilon_2} \cdots x_{i_{m-1}}^{\epsilon_{m-1}}
x_{i_m}^{(\epsilon_m +1)/2}) \in C_2^{\rm nh}(\pi_1(M) ;\Z). $$
Then, we can easily verify $ \partial_1^{\Delta} \circ c_1 = c_0 \circ \partial_1 $ and $ \partial_2^{\Delta} \circ c_2 = c_1 \circ \partial_2 $. Let $\mathcal{O}_M \in C_3( \widetilde{M};\Z) $ be the basis. Notice that $\partial^{\Delta}_2 \circ c_2 \circ \partial_3 (\mathcal{O}_M )= c_1 \circ \partial_2 \circ\partial_3 (\mathcal{O}_M ) =0$, that is, $c_2 \circ \partial_3 (\mathcal{O}_M ) $ is a 2-cycle. If we expand $c_2 \circ \partial_3 (\mathcal{O}_M ) $ as $\sum n_i (g_0^i, g_1^i, g_2^i)$ for some $ n_i \in \Z, g_j^i \in G $, then $\mathcal{O}_M ':=- \sum n_i (1, g_0^i, g_1^i, g_2^i)$ satisfies $\partial_3^{\Delta }(\mathcal{O}_M ') = c_2 \circ \partial_3 (\mathcal{O}_M ) $. Therefore, the correspondence $\mathcal{O}_M \mapsto \mathcal{O}_M ' $ gives rise to a chain map $c_3 : C_*(\widetilde{M} ) \ra C_*^{\rm Nor}(\pi_1(M); \Z ) $, as desired. In conclusion, the above discussion can be summarized as follows:
\begin{prop}\label{exh} For any homomorphism $f: \pi_1(M) \ra G $, the pushforward
$ f_* \circ \iota_*[M]$ is equal to $ 1 \otimes_{\pi_1(M)} f_* \circ c_3( \mathcal{O}_M) $ in $H_3^{\rm nh}(G;\Z) $.
\end{prop}
To conclude, if we know an explicit presentation of $\pi_1(M)$ and a representative of
the 3-cocycle $\psi:G^4 \ra A$, in principle, we can compute $ \DW_{\psi}(M) $ in terms of the chain map $c_*$ (with the help of computer program).

\subsection*{Acknowledgment}
The author expresses his gratitude to an anonymous referee
for reading this paper and giving him valuable comments.

\normalsize
\noindent
Department of Mathematics, Tokyo Institute of Technology 2-12-1 Ookayama, Meguro-ku Tokyo 152-8551 Japan


\begin{thebibliography}{99}

\small
\ifx\undefined\bysame
\newcommand{\bysame}{\leavevmode\hbox to3em{\hrulefill}\,}
\fi

\bibitem[Bro]{Bro} K. S. Brown, {\it Cohomology of Groups},
Graduate Texts in Mathematics, {\bf 87}, Springer-Verlag, New York, 1994.

\bibitem[BH]{BH} R. Brown, H. Huebschmann, {\it Identities among relations. Low-dimensional topology} (Bangor, 1979), pp. 153--202, London Math. Soc. Lecture Note Ser., 48, Cambridge Univ. Press, Cambridge-New York, 1982. MR0662431

\bibitem[BC]{BC}H. U. Boden, C. L. Curtis,
{\it Splicing and the $SL_2(\mathbb{C})$ Casson invariant}. Proc. Amer.
Math. Soc., 136 (2008): 2615--2623


\bibitem[CG]{CG}A. Casson, C. M. Gordon, {\it Cobordism of classical knots}, Progr. Math. 62, Birkh\"{a}user, Boston (1986) 181--199

\bibitem[CFH]{CFH} A. Conway,
S. Friedl, G. Herrmann, {\it Linking forms revisited}, Pure and Applied Mathematics Quarterly 12 (2016), 493--515
\bibitem[DW]{DW}
R. Dijkgraaf, E. Witten, {\it Topological gauge theories and group cohomology},
Comm. Math. Phys. {\bf 129} (1990) 393--429.

\bibitem[Kab]{Kab}Y. Kabaya,
{\it Cyclic branched coverings of knots and quandle homology}, Pacific Journal of Mathematics, {\bf 259} (2012), No. 2, 315--347

\bibitem[KKV]{KKV}
G. Kim, Y. Kim, A. Vesnin, {\it The knot $5_2$ and cyclically presented
groups}, J. Korean Math. Soc. {\bf 35} (1998), 961--980

\bibitem[Lea]{Lea}
I. J. Leary, {\it The mod-$p$ cohomology rings of some $p$-groups}. Math. Proc. Cambridge Philos. Soc. 112 (1992), no. 1, 63--75. MR1162933

\bibitem[Lic]{Lic} W.B. Lickorish, {\it An introduction to knot theory}, Springer-Verlag, Berlin - New York, 1974.
\bibitem[Lyn]{Lyn} R. Lyndon, {\it Cohomology theory of groups with a single defining relation}, Ann. of Math. {\bf 52} (1950), 650--665
\bibitem[Mar]{Mar}
S. T. Martins, {\it Diagonal approximation and the cohomology ring of torus fiber bundles},
Internat. J. Algebra Comput. 25 (2015), no. 3, 493--530


\bibitem[MP]{MP}
A. Miller and M. Powell, {\it Symmetric chain complexes, twisted Blanchfield pairings
and knot concordance}. Algebr. Geom. Topol., 18(2018), 3425--3476

\bibitem[No1]{Nos} T. Nosaka, {\it An $SL_2(\R) $-Casson invariant and Reidemeister torsions}, preprint, available at arXiv:2005.06132
\bibitem[No2]{No2}
T. Nosaka, {\it Quandle cocycle invariants of knots using Mochizuki's 3-cocycles and Dijkgraaf-Wittten invariants of 3-manifolds},
Algebraic and Geometric Topology 14 (2014) 2655--2692

\bibitem[Sc]{Sc} P. Scott, {\it There are no fake Seifert fibre spaces with infinite $\pi_1$}, Ann. of Math. {\bf 117} (1983), 35--70

\bibitem[Sei]{Sei} H. Seifert, {\it Verschlingungs invarianten}, Sitzungsber. Preu\ss. Akad. Wiss., Phys.-Math. Kl. 1933, No.26-29, (1933) 811--828

\bibitem[Sie]{Sie1} A. J. Sieradski, {\it Combinatorial squashings, 3-manifolds, and the third homology of groups}, Invent. math. {\bf 84}(1986), 121--139


\bibitem[Tro]{Tro} H. F. Trotter,
{\it Homology of group systems with applications to knot theory}, Ann. of Math. {\bf 76}
(1962), 464--498

\bibitem[Wall]{Wall}
T. C. Wall. {\it Classification problems in differential topology. VI. Classification
of $(s-1)$-connected $(2s + 1)$-manifolds}. Topology, {\bf 6} (1967): 273--296

\bibitem[Wakui]{W2}M. Wakui, {\it On Dijkgraaf-Witten invariant for $3$-manifolds}. Osaka J. Math. 29 (1992), no. 4, 675--696. MR1192735

\bibitem[Waki]{W} N. Wakijo, {\it Twisted Reidemeister torsions of 3-manifolds via Heegaard splittings}, Topology Appl. 299 (2021), 107731. MR4270619

\bibitem[Wise]{BW}
D. T. Wise, {\it The structure of groups with a quasiconvex hierarchy}, Electronic Research Announcements In Mathematical Sciences, V. 16, 44--55

\end{thebibliography}
\end{document}